\documentclass{ifacconf}

\usepackage{graphicx}      
\usepackage{natbib}        

\usepackage{amsmath,amssymb,amsfonts}
\usepackage{algorithmic}
\usepackage{subfigure,graphicx}
\usepackage{textcomp}
\usepackage{tikz}

\newtheorem{theorem}{Theorem}[section]
\newtheorem{corollary}[thm]{Corollary}
\newtheorem{lemma}[thm]{Lemma}
\newtheorem{remark}[thm]{Remark}
\newtheorem{assumption}[theorem]{Assumption}    
\newcommand{\Sym}{\mathrm{Sym}}
\newcommand{\Skew}{\mathrm{Skew}}
\newcommand{\SO}{\operatorname{SO}(3)}

\newcommand{\meas} {{F}}

\begin{document}
\begin{frontmatter}

\title{A Global, Continuous, and Exponentially Convergent Observer for Gyro Bias and Attitude of a Rigid Body\thanksref{footnoteinfo}} 

\thanks[footnoteinfo]{The first author has been supported in part by the KUSTAR-KAIST Institute, KAIST, Korea,  by KAIST under grants G04170001 and N11180231.
The second author has been supported in part by NSF under the grants 1837382 and 1760928, and by AFOSR under the grant FA9550-18-1-0288.}

\author[First]{Dong Eui Chang} 
\author[Second]{Taeyoung Lee}

\address[First]{School of Electrical Engineering,  Korea Advanced Institute of Science and Technology, Daejeon, 34141, Korea
         (e-mail: dechang@kaist.ac.kr)}
\address[Second]{Department of Mechanical and Aerospace Engineering,  George Washington University, Washington DC, 20052, USA (e-mail: tylee@gwu.edu)}

\begin{abstract}                
We propose a 12-dimensional, global, continuous, and exponentially convergent observer for  gyro bias and attitude of a rigid body. Any attitude observer developed on the special orthogonal group suffers from the topological restriction that prohibits global attractivity in continuous flow. In this paper, the observer is designed in the set of 3 by 3 real matrices, thus making the topological obstruction on the special orthogonal group irrelevant. The efficacy of the proposed approach against other attitude observers is illustrated by an indoor experiment utilizing visual landmarks.  
\end{abstract}

\begin{keyword}
Gyro bias, estimator, observer, global convergence, rigid body
\end{keyword}

\end{frontmatter}

\section{INTRODUCTION}

Estimating the attitude of a rigid body from vector measurements has been for decades a problem of interest, because of its importance for a variety of technological applications such as satellites or unmanned aerial vehicles.  In many cases, the measurement of angular velocity  is corrupted with a sensor bias, so it is important to estimate the bias in the angular velocity measurement. This issue can  be more precisely expressed in mathematical terms. 
Consider the rigid body kinematics equation
\begin{equation}\label{rbk}
\dot R = R\hat \Omega,
\end{equation}
where $R \in \SO$ is the $3\times 3$ rotation matrix and $\Omega \in \mathbb R^3$ is the body angular velocity of the rigid body. Here, the hat map $\hat{} : \mathbb R^3 \rightarrow \mathfrak{so}(3)$ is defined in \eqref{def:Omega} in Appendix.   Assume that we measure the angular velocity  with a sensor that contains a bias $b$ as follows:
\[
\Omega_{\rm m} = \Omega + b
\]
where $\Omega_{\rm m}$ is the measured angular velocity and $b\in \mathbb R^3$ is the bias vector. It is often the case that the bias $b$ is slowly varying in time, so it is common to assume that it is constant, which can be expressed as
\begin{equation}\label{constant:b}
\dot b = 0.
\end{equation}
The task is to estimate $(R,b)$ with measurements of vector quantities such as gravity, geomagnetic field, landmarks, etc. Estimating $R$ is relatively straightforward whereas estimating $b$ is a challenge.

In this paper, we will build a nonlinear observer for the system \eqref{rbk} and \eqref{constant:b} to estimate $(R,b)$ from vector measurements. A noticeable technique employed here is an embedding technique such that we first embed the rigid body kinematics into $\mathbb R^{3\times 3}$ since $\SO$ is an embedded submanifold of $\mathbb R^{3\times 3}$. We then build an observer in the extended state space $\mathbb R^{3\times 3} \times \mathbb R^3$ instead of $\SO \times \mathbb R^3$ to estimate $(R,b)$. This embedding technique allows us to get around the notorious topological defect of $\SO$, which  hinders global convergence on $\SO$, such that we achieve a globally exponentially convergent observer for $(R,b)$. The embedding technique has been developed and applied successfully in various areas of control and numerics; optimal control \citep{Ch11},  tracking controller design \citep{Ch18,Ch18b},  and structure-preserving integration\citep{ChJiPe16}.

The design of angular velocity bias observer for a rigid body has been  studied by many researchers  \citep{GrFoJoSa12, Ma08, MaSa07, MaSa10, WuKaLe15, BaSiOl12a, BaSiOl12b,BaSiOl12c,GrFoJoSa15,MaSa17,BeAbTa17}.  In \cite{GrFoJoSa12}, \cite{Ma08}, \cite{MaSa07}, \cite{MaSa10}, \cite{WuKaLe15}, and \cite{BeAbTa17},  observers are built on $\SO$ or unit quaternions, where the seminal work by \cite{Ma08} introduces the formal observers on $\SO$ with theoretical guarantees for the first time. Due to the incontractibility of $\SO$ and the set of unit quaternions, none of the continuous observers  are globally convergent whereas the hybrid observer with switching rules in \cite{BeAbTa17}, which gets around the topological obstruction on $\SO$, is globally convergent.  

Alternatively, there is another family of observers which are built in Euclidean space \citep{BaSiOl12a,BaSiOl12b,BaSiOl12c,GrFoJoSa15,MaSa17} so that the topological obstruction on $\SO$ can be bypassed, among which  \cite{BaSiOl12a} and \cite{BaSiOl12b} are  the seminal works introducing the approach using Euclidean space rather than $\SO$ or unit quaternions.   In \cite{BaSiOl12c}, a cascade observer with global exponential stability is proposed, where several desirable properties, such as a complementary structure are achieved.  However, the dimension of the cascade observer  increases as the more vector measurements become available. 
More specifically, the dimension of the observer is $3N+12$ when there are measurements of $N$ vectors. However, the redundancy can be interpreted positively.  In contrast to this,  the observer presented in \cite{MaSa17} is economically designed with only two distinct vectors so that the dimension of the observer is only 9, and it is proven to be uniformly globally asymptotically convergent  and locally exponentially convergent, but not globally exponentially convergent. A globally exponentially stable attitude observer is proposed in~\cite{GrFoJoSa15} with application to GNSS/INS integration. 
While it is claimed to be globally exponentially convergent, this observer requires  knowledge of an upper bound of the magnitude of  unknown constant gyro bias that is used in such a way that the bound needs to be reset to accommodate   biases that are bigger than the initially chosen upper bound. Hence, at a closer look it is not hard to see that the convergence of their observer is only semi-global.

%

We here propose a continuous and globally exponentially convergent observer for attitude and  gyro bias for a rigid body that is constructed on $\mathbb R^{3\times 3} \times \mathbb R^3$ instead of $\SO \times \mathbb R^3$ to avoid the topological restriction on $\SO$. However, compared with the aforementioned attitude observers in  \cite{BaSiOl12a}, \cite{BaSiOl12b}, \cite{BaSiOl12c},  \cite{GrFoJoSa15}, and \cite{MaSa17}, there are unique, desirable features in our observer. First,  our approach respects the matrix operations on $\mathbb R^{3\times 3}$.  Second, the dimension of our observer is always 12 irrespective of the number of vector measurements. Further, it does not require any knowledge of an upper bound of the magnitude of  bias, and it does not involve any hybrid switching rule.  Finally, it exhibits global exponential convergence for the attitude and the gyro bias estimation errors. 

The observer in this paper is first proposed in a unified form, and then various specific forms of observer are derived from it, including the case when the reference directions are time-varying.    
Our proof for global and exponential convergence of the observer is straightforward and easy to understand, which does not require any persistent excitation assumption. 
Experimental results are also presented to illustrate an excellent performance of our global observer in comparison with the observer in the seminal work by \cite{Ma08}, where the comparison is fair since the two observers have almost the same structure so that the same gain values can be used.

\section{MAIN RESULTS}
We first invite the reader to read the Appendix  to get acquainted with  the mathematical preliminaries  that will be used throughout the paper. 
The kinematic equation of a rigid body is given by
\begin{equation}\label{rigid:body}
\dot R = R \hat \Omega,
\end{equation}
where $R \in \SO$ is the rotation or attitude of a rigid body,  $\Omega \in \mathbb R^3$ is the body angular velocity, and the symbol $\wedge$ over $\Omega$ denotes the hat map, $\wedge : \mathbb R^3 \rightarrow \mathfrak{so}(3)$,  defined in the Appendix.   We make the following three assumptions.
\begin{assumption}\label{assumption:A}
A matrix-valued signal $A \in \mathbb R^{3\times 3}$ is available and  can be expressed as
\begin{equation}\label{relation:AG}
A = \meas R,
\end{equation}
 where $\meas$ is a constant  invertible matrix in $\mathbb R^{3\times 3}$ and $R \in \SO$ is the attitude of the rigid body.
\end{assumption}

\begin{assumption}\label{assumption:Omega}
A measured body angular velocity $\Omega_{\rm m}$ with bias is available and related to  the body angular velocity $\Omega$ of the rigid body as follows:
\begin{align*}
 \Omega_{\rm m} = \Omega + b,
\end{align*}
where $b$ is an unknown bias vector.  
\end{assumption}

\begin{assumption}\label{assumption:Ombbound}
The trajectory of angular velocity $\Omega(t)$ is bounded, and the bias vector $b$ is constant. 
\end{assumption}

We propose the following  observer:
\begin{subequations}\label{filter:Chang}
\begin{align}
\dot{\bar A} &= \bar A {\hat \Omega}_{\rm m} - A \hat{\bar b} +k_P  (A - \bar A),\\
\dot{\bar b} &= k_I \Skew (A^T \bar A)^\vee \label{filter:Chang:b}
\end{align}
\end{subequations}
with  $k_P>0$ and $k_I >0$, where  $(\bar A, \bar b) \in \mathbb R^{3\times 3} \times \mathbb R^3$ is an estimate of $(A,b)$. So,  $(\meas^{-1}\bar A, \bar b)$ becomes an estimate of $(R,b)$ by Assumption \ref{assumption:A}. The global and exponentially convergent property of this observer is proven in the following theorem. 

\begin{theorem} \label{main:theorem:filter}
Let
\[
E_A = A-\bar A, \quad e_b = b - \bar b.
\]
Under Assumptions \ref{assumption:A} -- \ref{assumption:Ombbound}, for any $k_P>0$ and $k_I>0$
there exist numbers $a>0$ and $C>0$ such that
\begin{equation}\label{final:eqn}
\| E_A(t)\| + \|e_b(t)\| \leq C( \| E_A(0)\| + \|e_b(0)\| )e^{-at}
\end{equation}
for all $t\geq0$ and all $(\bar A(0), \bar b(0)) \in \mathbb R^{3\times 3} \times \mathbb R^3$. 

\begin{pf}
From \eqref{rigid:body} and Assumption \ref{assumption:A}, $A(t)$ satisfies
\begin{equation}\label{Adot:eq}
\dot A =  A \hat \Omega. 
\end{equation}
By Assumption \ref{assumption:Omega}, the observer can be written as
\begin{subequations}\label{modified:filter:Chang}
\begin{align}
\dot{\bar A} &= \bar A(\hat \Omega + \hat b) - A \hat{\bar b} +k_P E_A,\\
\dot{\bar b} &= -k_I \Skew (A^T E_A)^\vee
\end{align}
\end{subequations}
since $\Skew (A^T \bar A)= -\Skew(A^TE_A)$.  By Assumption \ref{assumption:Ombbound}, there are numbers $B_\Omega>0$ and $B_b>0$ such that $\|\Omega (t)\| \leq B_\Omega$ for all $t\geq 0$ and $\|b\|\leq B_b$. Let $B = \max\{B_\Omega, B_b\}$. Then, there is a number $\epsilon$ such that
\[
0 < \epsilon <\frac{1}{\|\meas\| \sqrt{k_I}}
\]
and
\[
0 < \epsilon <  \frac{4k_P \lambda_{\rm min}(\meas^T\meas)}{ \|\meas\|^2(4k_I  \lambda_{\rm min}(\meas^T\meas) +  (k_P + 3\sqrt{2} B)^2 )},
\]
where $ \lambda_{\rm min}(\meas^T\meas)$ denotes the smallest  eigenvalue of $\meas^T\meas$, which is positive since $\meas$ is invertible. 
The following three quadratic functions of $(\|E_A\|, \|e_b\|)$ are then all positive definite:
\begin{align*}
V_1(\|E_A\|, \|e_b\|) &= \frac{1}{2} \|E_A\|^2 + \frac{1}{k_I}\|e_b\|^2 \! - \! \sqrt{2}\epsilon \|\meas\| \|E_A\| \|e_b\|,\\
V_2(\|E_A\|, \|e_b\|) &= \frac{1}{2} \|E_A\|^2 + \frac{1}{k_I}\|e_b\|^2 \! + \! \sqrt{2}\epsilon \|\meas\| \|E_A\| \|e_b\|,\\
V_3(\|E_A\|, \|e_b\|) &= (k_P - \epsilon k_I \|\meas\|^2) \|E_A\|^2 \\
&\quad+ 2\epsilon  \lambda_{\rm min}(\meas^T\meas)\|e_b\|^2 \\
&\quad -\epsilon (\sqrt{2}k_P + 6B)\|\meas\|\|E_A\| \|e_b\|.
\end{align*}
Hence, there are numbers $\alpha>0$ and $\beta >0$ such that 
\begin{equation}\label{V213}
V_2 \leq \alpha V_1, \quad \beta V_2 \leq  V_3.
\end{equation}
Let 
\[
V(E_A, e_b) = \frac{1}{2} \|E_A\|^2 + \frac{1}{k_I} \|e_b\|^2 + \epsilon \langle E_A, A\hat e_b \rangle,
\]
which satisfies
\begin{equation}\label{V1:V:V2}
V_1(\|E_A\|, \|e_b\|) \leq V(E_A, e_b) \leq V_2(\|E_A\|, \|e_b\|)
\end{equation}
for all $(E_A, e_b) \in \mathbb R^{3\times 3} \times \mathbb R^3$ by the Cauchy-Schwarz inequality,  statements 3 and 5 in Lemma \ref{lemma:prelim} in the Appendix, and $\|A\| = \|\meas R\| = \|\meas\|$ since $R\in \SO$. From  \eqref{Adot:eq}, \eqref{modified:filter:Chang},  and the assumption of the bias $b$ being constant, it follows that the estimation error $(E_A, e_b)$ obeys
\begin{align*}
\dot E_A &= E_A(\hat \Omega + \hat b) - A\hat e_b - k_PE_A,\\
\dot e_b &= k_I \Skew(A^TE_A)^\vee.
\end{align*}
Along any trajectory of the composite system consisting of the rigid body \eqref{rigid:body} and the observer \eqref{filter:Chang},
\begin{align*}
\frac{dV}{dt} &= \langle E_A, E_A(\hat \Omega + \hat b) - A\hat e_b - k_PE_A\rangle \\
&\quad+ 2\langle  e_b, \Skew (A^TE_A)^\vee \rangle \\
&\quad + \epsilon \langle E_A (\hat \Omega + \hat b) - A\hat e_b - k_PE_A, A\hat e_b\rangle  \\
&\quad+ \epsilon \langle E_A, A\hat \Omega \hat e_b\rangle +\epsilon k_I \langle E_A, A\Skew(A^TE_A)\rangle \\
&\leq -(k_P - \epsilon k_I \|\meas\|^2) \|E_A\|^2 - 2\epsilon  \lambda_{\rm min}(\meas^T\meas)\|e_b\|^2 \\
&\quad +\epsilon (\sqrt{2}k_P + 6B)\|\meas\|\|E_A\| \|e_b\| \\
&= -V_3 \leq -\beta V_2 \leq - \beta V,
\end{align*}
where the following have been used:
\begin{align*}
&\langle E_A, E_A(\hat \Omega + \hat b) \rangle =   \langle E_A^TE_A, (\hat \Omega + \hat b) \rangle =0,\\
&  \langle E_A, A\hat e_b\rangle =  \langle \Skew(A^TE_A), \hat e_b\rangle = 2\langle \Skew (A^TE_A)^\vee, e_b\rangle,\\
&\langle A\hat e_b, A\hat e_b\rangle \geq \lambda_{\rm min}(\meas^T\meas)\|R\hat e_b\|^2 = 2\lambda_{\rm min}(\meas^T\meas)\| e_b\|^2,\\
&\langle E_A, A\Skew(A^TE_A)\rangle = \|\Skew(A^TE_A)\|^2\\
&\qquad\qquad\qquad\qquad \quad \leq \|A^TE_A\|^2 \leq \|\meas\|^2 \|E_A\|^2.
\end{align*}
Hence, $V(t) \leq V(0) e^{-\beta t}$ for all $t\geq 0$ and all $(\bar A(0), \bar b(0)) \in \mathbb R^{3\times 3} \times \mathbb R^3$. It follows from \eqref{V213} and \eqref{V1:V:V2} that
\[
V_1(t) \leq V(t) \leq V(0) e^{-\beta t}  \leq V_2(0) e^{-\beta t}  \leq \alpha V_1(0) e^{-\beta t} 
\]
for all $t\geq 0$ and all $(\bar A(0), \bar b(0)) \in \mathbb R^{3\times 3} \times \mathbb R^3$. Since $0<\epsilon < 1/ (\|\meas\|\sqrt{k_I})$,  the map defined by
\[
(x_1, x_2) \mapsto \sqrt{\frac{1}{2}x_1^2 + \frac{1}{k_I}x_2^2 - \sqrt 2 \epsilon \|\meas\| x_1x_2}
\]
is a norm on $\mathbb R^2$, where $(x_1,x_2) \in \mathbb R^2$, which is equivalent to the 1-norm on $\mathbb R^2$ since all norms are equivalent on a finite-dimensional vector space. Hence,  $V_1(t) \leq \alpha V_1(0)e^{-\beta t}$ implies that there exists $C>0$ such that \eqref{final:eqn} holds  for all $t\geq 0$ and all $(\bar A(0), \bar b(0)) \in \mathbb R^{3\times 3} \times \mathbb R^3$, where $a = \beta/2$. 
\end{pf}
\end{theorem}

 \begin{remark}
 1. Notice in the proof of Theorem \ref{main:theorem:filter} that the numbers $a$ and  $C$  in \eqref{final:eqn} may depend on $B_\Omega$ and $B_b$, but it has not prevented us from showing the exponential convergence of the observer. Moreover, the choice of $k_P$ and $k_I$ is totally independent of $B_\Omega$ and $B_b$.
 
 2. One can easily generalize the form of observer in \eqref{filter:Chang} by replacing the scalar gains $k_P$ and $k_I$ with matrix gains. 
  \end{remark}

\begin{corollary} \label{main:corollary:filter}
Suppose that Assumptions \ref{assumption:A} -- \ref{assumption:Ombbound} hold, and  let
\[
E_R = R-\meas^{-1}\bar A, \quad e_b = b - \bar b.
\]
 Then, 
there exist numbers $a>0$ and $C>0$ such that
\begin{equation}\label{final:eqn:2}
\| E_R(t)\| + \|e_b(t)\| \leq C( \| E_R(0)\| + \|e_b(0)\| )e^{-at}
\end{equation}
for all $t\geq0$ and all $(\bar A(0), \bar b(0)) \in \mathbb R^{3\times 3} \times \mathbb R^3$.
\begin{pf}
Use $\|E_R\| / \|\meas^{-1}\| \leq\|E_A\| \leq \|\meas\| \|E_R\|$ and  \eqref{final:eqn} with the constant $C$ redefined appropriately.
\end{pf}
\end{corollary}
In other words, the estimate $(\meas^{-1}\bar A(t), \bar b(t))$ of the pair of attitude and gyro bias converges  globally and exponentially  to the true value  $(R(t), b)$ as $t$ tends to infinity. 

\begin{remark}
If $\meas=I$, then the observer \eqref{filter:Chang} reduces to 
\begin{subequations}\label{filter:Chang:R}
\begin{align}
\dot{\bar R} &= \bar R {\hat \Omega}_{\rm m} - R \hat{\bar b} +k_P  (R - \bar R),\\
\dot{\bar b} &= k_I \Skew (R^T \bar R)^\vee, 
\end{align}
\end{subequations}
 where  $(\bar R, \bar b) \in \mathbb R^{3\times 3} \times \mathbb R^3$ is an estimate of $(R,b)$.   The global  observer in \eqref{filter:Chang:R} may look similar to but is different from the non-global  observers 
proposed by Mahony  et al that appear in   (12) and (13) in \cite{Ma08}.
 \end{remark}

 \begin{remark}
 We can relax Assumption \ref{assumption:A} by allowing the matrix $\meas$ to be time-varying. More specifically, we make the following assumption: there  are numbers $\ell_{\rm min} >0$  and $\ell_{\rm max} >0$ such that 
 \begin{equation}\label{relax:G}
 \ell_{\rm min} \leq  \lambda_{\rm min} (\meas^T(t)\meas(t)) \leq \lambda_{\rm max} (\meas^T(t)\meas(t)) \leq \ell_{\rm max}
 \end{equation}
  for all $t \geq 0$. 
In this case, we propose the following  observer:
\begin{align*}
\dot{\bar A} &= \bar A {\hat \Omega}_{\rm m} - A \hat{\bar b} +k_P  (A - \bar A) + \dot \meas \meas^{-1}A,\\
\dot{\bar b} &= k_I \Skew (A^T \bar A)^\vee
\end{align*}
with  $k_P>0$ and $k_I >0$, where  $(\bar A, \bar b) \in \mathbb R^{3\times 3} \times \mathbb R^3$ is an estimate of $(A,b)$.  It is not difficult to show that Theorem  \ref{main:theorem:filter} and Corollary \ref{main:corollary:filter} still hold for this observer with the relaxed assumption on $\meas(t)$ as above.  The knowledge on the values of $\ell_{\rm min} $  and $\ell_{\rm max}$ is not required here. The proof would involve a small modification of the proof of Theorem  \ref{main:theorem:filter}, which is left to the reader. 
 \end{remark}

\begin{remark}
 Instead of \eqref{filter:Chang}, we can  consider the following form of observer:
 \begin{subequations}\label{filter:Chang:inverse}
\begin{align}
\dot{\bar A} &= \bar A {\hat \Omega}_{\rm m} - A \hat{\bar b} +k_P  (A - \bar A),\label{filter:Chang:inverse:a}\\
\dot{\bar b} &= k_I \Skew (A^{-1} \bar A)^\vee
\end{align}
\end{subequations}
with  $k_P>0$ and $k_I >0$, where the only difference between \eqref{filter:Chang} and \eqref{filter:Chang:inverse} is in the  equation for $\dot {\bar b}$.  It is not difficult to prove that this observer also enjoys the property of global and exponential convergence for any $k_P>0$ and $k_I >0$, whose proof is left to the reader.
 \end{remark}

 We now derive from \eqref{filter:Chang} various  observers of {\it concrete} form that estimate $(R,b)$ from vector measurements. 
  Assume that there is a set $\mathcal S = \{s_i, 1\leq i \leq m\}$ of $m$ known fixed inertial vectors, where each $s_i$ in $\mathcal S$ is a vector in $\mathbb R^3$, such that the rank of $\mathcal S$ is 3. If the rank of $S$ is only 2, then pick any two mutually independent vectors $s_i$ and $s_j$ from $\mathcal S$ and add $s_i \times s_j$ to $\mathcal S$, which will make $S$ have rank 3.  Assume also that measurements of  vectors are made in the body-fixed frame and the set of the measured vectors is denoted by $\mathcal C = \{ c_i , 1\leq i \leq m\}$ and  related to $\mathcal S$ as follows:
 \[
 c_i = R^Ts_i, \quad i = 1, \ldots, m,
 \]
 where $R$ is the attitude of the rigid body.  Let
 \begin{equation}\label{def:SC}
 S =  \begin{bmatrix}
 s_1& \cdots & s_m
 \end{bmatrix}, \quad C =  \begin{bmatrix}
 c_1& \cdots & c_m
 \end{bmatrix}
 \end{equation}
  be $3 \times m$ matrices made of the {\it column} vectors from $\mathcal S$ and ${\mathcal C}$.

 \begin{corollary}\label{corollary:obs:linear}
Let  $S$ and $C$ be given  in \eqref{def:SC}, and let 
\[
W = \begin{bmatrix}
w_1 & \cdots & w_m
\end{bmatrix} \in \mathbb R^{3\times m},
\]
where $w_i \in \mathbb R^3$ denotes the $i$th column vector of $W$ for $1\leq i \leq m$. Suppose
 \begin{equation}\label{linear:GA}
 \meas = WS^T, \quad A = WC^T
 \end{equation}
 such that the rank of $\meas$ is 3. Then, \eqref{relation:AG} is satisfied, and \eqref{filter:Chang} reduces to 
\begin{subequations}\label{obs:linear}
\begin{align}
\dot{\bar A} &= \bar A {\hat \Omega}_{\rm m} - A \hat{\bar b} +k_P  (A - \bar A),\label{obs:linear:A}\\
 \dot{\bar b}  &=- k_I \sum_{i=1}^m c_i \times \bar A^T w_i,  \label{obs:linear:b}
 \end{align}
 \end{subequations}
 where $k_I/2$ has been replaced with $k_I$ in \eqref{obs:linear:b} to make \eqref{obs:linear:b} look simple.
 \begin{pf}
Trivial by statement 6 of Lemma \ref{lemma:prelim} in the Appendix.
 \end{pf}
 \end{corollary}

 \begin{corollary}\label{corollary:obs:quad}
 Suppose 
 \begin{equation}\label{quad:GA}
 \meas = SWS^T, \quad   A = SWC^T,
 \end{equation}
 where $W = [w_{ij}]$ is an $m\times m$ matrix  such that $\meas$ has rank 3. Then, \eqref{relation:AG} is satisfied, and \eqref{filter:Chang} reduces to 
\begin{subequations}\label{obs:quad}
\begin{align}
\dot{\bar A} &= \bar A {\hat \Omega}_{\rm m} - A \hat{\bar b} +k_P  (A - \bar A),\label{obs:quad:A}\\
 \dot{\bar b}  &=- k_I \sum_{i=1}^m \sum_{j=1}^m w_{ij} c_j \times \bar A^T s_i,  \label{obs:quad:b}
 \end{align}
 \end{subequations}
 where $k_I/2$ has been replaced with $k_I$ in \eqref{obs:quad:b} to make \eqref{obs:quad:b} look simple.

 \end{corollary}
  
 \begin{remark}
 If $W$ is a diagonal matrix in Corollary \ref{corollary:obs:quad}, then  $\meas$ and $A$ in \eqref{quad:GA} become $
 \meas = \sum_{i=1}^m w_{ii} s_is_i^T$ and $A = \sum_{i=1}^m w_{ii} s_i c_i^T$,
 and the observer  \eqref{obs:quad} reduces to 
 \begin{subequations}\label{obs:quad:var}
\begin{align}
\dot{\bar A} &= \bar A {\hat \Omega}_{\rm m} - A \hat{\bar b} +k_P  (A - \bar A),\label{obs:quad:var:A}\\
 \dot{\bar b}  &=- k_I \sum_{i=1}^m w_{ii} c_i \times \bar A^T s_i.  \label{obs:quad:var:b}
 \end{align}
 \end{subequations}
The observer \eqref{obs:quad:var}, which is global, may look similar to   but is different from the non-global observer in (32) in \cite{Ma08}.  A study by simulation and experiment later in this paper will demonstrate a noticeable difference in performance between the two observers.
 \end{remark}

\begin{remark}
 Corollary \ref{corollary:obs:quad} can be regarded as a special case of  Corollary \ref{corollary:obs:linear} since the substitution of $SW$ into $W$ in \eqref{linear:GA} would yield \eqref{quad:GA}.  Likewise, if we choose $W$ in Corollary  \ref{corollary:obs:linear} such that it depends nonlinearly on $S$, then it will produce another observer that is different from \eqref{obs:linear} and \eqref{obs:quad}.  Also, if we allow $W$ in Corollary \ref{corollary:obs:linear} and Corollary \ref{corollary:obs:quad} to vary in time, say by making it depend on $c_i(t)$'s or $A(t)$ or  allowing $s_i$ to vary in time, such that \eqref{relax:G} is satisfied for all $t\geq 0$,  then $\dot{\bar A}$ dynamics in \eqref{obs:linear:A} and \eqref{obs:quad:A} modifies to 
\[
\dot{\bar A} = \bar A {\hat \Omega}_{\rm m} - A \hat{\bar b} +k_P  (A - \bar A) + \dot \meas \meas^{-1}A
\]
while \eqref{obs:linear:b} and \eqref{obs:quad:b} remain unchanged.
\end{remark}
  
 \begin{remark}
 Our observer, whose dynamics evolve globally in Euclidean space, is straightforward to numerically integrate, whereas most observers on $\SO$  would require an  operation of projection onto $\SO$  at each step of  numerical integration, which  adds numerical errors to integration.  The use of unit quaternions for numerical integration, which also requires  projections from $\mathbb R^4$ onto the set of unit quaternions ${\rm S}^3$, is not immune to such possibly accumulative numerical errors in integration. Notice that just projecting points from $\mathbb R^4$ onto ${\rm S}^3$ does not guarantee precise projections onto the correct points in the trajectory on ${\rm S}^3$. It only guarantees that   the projected image is in ${\rm S}^3$, and this kind of errors may accumulate during integration. In contrast, our observer, which runs in Euclidean space, does not have any such problem at all in numerical integration.  The same merit applies to the observers in \cite{BaSiOl12a,BaSiOl12b,BaSiOl12c,GrFoJoSa15,MaSa17}.
  \end{remark}
 
 \begin{remark}
 Putting the numerical integration issue aside, we can approximate the trajectory of estimates $\bar R(t) \in \mathbb R^{3\times 3}$ with a curve of rotation matrices using polar decomposition.   Since the operation of polar decomposition is continuous \citep[Theorem 4.1.4]{Fa08}, the $\SO$ part of $\bar R(t)$ obtained from polar decomposition also converges to $R(t)$ as $t\rightarrow \infty$. Any error due to this approximation is not accumulative during integration because the approximate value is not used in the numerical integration of observer dynamics but is only fed to controllers. 
 However, this approximation may not even be  necessary when $\bar R(t)$ is directly used in feedback control.  Suppose that we have  an $\mathbb R^3$-valued control law $u(R,\Omega)$ for a rigid body system, where $R \in \SO$ is the attitude of the rigid body and $\Omega \in \mathbb R^3$ the body angular velocity. We can naturally extend the function $u(R,\Omega)$  to $\mathbb R^{3\times 3} \times \mathbb R^3$ by treating $R$ as a $3\times 3$ matrix after replacement of any occurrence of $R^{-1}$ with $R^T$ in the expression of $u(R,\Omega)$; refer to \cite{Ch18} for this extension. Recall that $(\bar R (t), \bar \Omega (t))$ with $\bar \Omega (t) = \Omega_{\rm m}(t) - \bar b(t)$ converges exponentially to $(R(t), \Omega (t))$. Hence,   if $(\bar R (t), \bar \Omega (t)) \approx (R(t), \Omega (t))$ in $\mathbb R^{3\times 3} \times \mathbb R^3$, then $u(R(t), \Omega (t)) \approx u(\bar R(t), \bar \Omega (t))$ in $\mathbb R^3$ by continuity of $u (\cdot , \cdot )$.  
 \end{remark}

\section{EXPERIMENTAL RESULTS}

We  validate the proposed observer  \eqref{filter:Chang} with experiments, and compare it with the Mahony observer that is given  in (32) in \cite{Ma08}. The reason why we make comparison with the Mahony observer is that it almost has the same structure as the proposed observer  \eqref{filter:Chang} such that fare comparison can be made with the same set of gains. Comparison with any other observer would raise an issue of fairness.

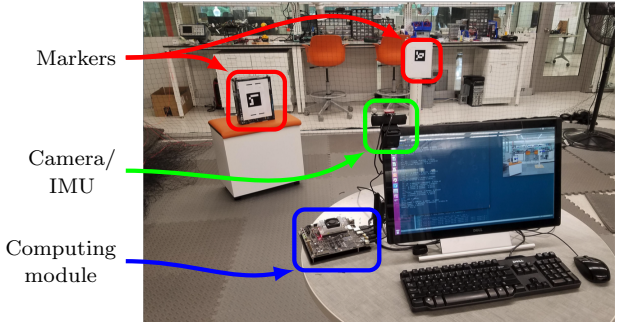
\begin{figure}[t]
\begin{tikzpicture}
\node [anchor=east] (marker) at (0.3,3.5) {\footnotesize Markers};
\node [anchor=east] (imu) at (0.3,2) {\footnotesize \shortstack[c]{Camera/\\IMU}};
\node [anchor=east] (jetson) at (0.3,0.8) {\footnotesize \shortstack[c]{Computing\\module}};
\begin{scope}[xshift=0.06\columnwidth]
    \node[anchor=south west,inner sep=0] (image) at (0,0) {\includegraphics[width=0.7\columnwidth]{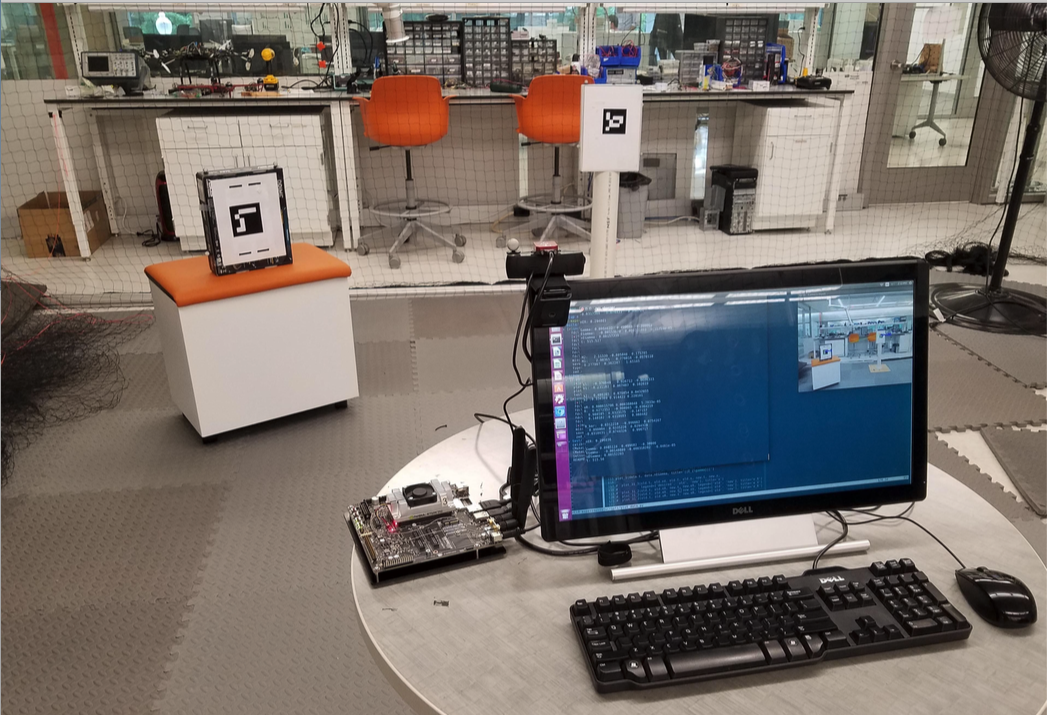}};
    \begin{scope}[x={(image.south east)},y={(image.north west)}]
        \draw[red,ultra thick,rounded corners] (0.18,0.6) rectangle (0.3,0.78);
        \draw[red,ultra thick,rounded corners] (0.55,0.75) rectangle (0.63,0.88);
        \draw [-latex, ultra thick, red] (marker) to[out=0, in=150] (0.18,0.78);
        \draw [-latex, ultra thick, red] (marker) to[out=0, in=150] (0.55,0.88);
        \draw[green,ultra thick,rounded corners] (0.47,0.53) rectangle (0.58,0.69);
        \draw [-latex, ultra thick, green] (imu) to[out=0, in=210] (0.47,0.53);
        \draw[blue, ultra thick,rounded corners] (0.32,0.16) rectangle (0.5,0.35);
        \draw [-latex, ultra thick, blue] (jetson) to[out=0, in=190] (0.32,0.16);
    \end{scope}
\end{scope}
\end{tikzpicture}
    \caption{Hardware configuration for experiment}\label{fig:exp_setup}
\end{figure}

The hardware system utilized in this paper is composed of the following parts: computing module (NVIDIA Jetson TX2) for sensor acquisition, observer implementation, and data logging; IMU (VectorNav VN-100) for angular velocity and acceleration measurements; camera (Logitech C930e) for line of sight measurements toward feature points; motion capture (VICON) system; see Figure \ref{fig:exp_setup}.
The camera and the IMU are rigidly attached with each other using a double-sided tape. 

We test our observer given in \eqref{obs:quad:var} with $m=3$ measurements.
Two distinct, square markers are placed in the field of view of the camera so that each marker is detected by using the OpenCV library \citep{GarMun14}. 
Then, the location of the markers in the image plane is converted into the corresponding unit-vector in the body-fixed frame.
The third measurement is the acceleration vector measured from IMU, which is considered as the direction of gravity. 


To obtain the corresponding directions in the inertial frame, the locations of the square markers and the camera/IMU are measured from the external motion capture system. 
The attitude of the combined camera/IMU body is also measured from the motion capture system, and it is considered as the true attitude against which the estimated attitude is compared. 
Also, to test the effects of a relatively large gyro bias over a short-time period, an artificial bias $b=(0,0.1,-0.2)$ is added to the angular velocity measurement from IMU. 

The proposed observer and the Mahony observer are implemented in C++ with multiple threads executing data acquisition, image processing, observer update, and data logging simultaneously at $50\,\mathrm{Hz}$. 

\begin{figure}
    \centerline{
        \subfigure[Attitude estimation error]{\includegraphics[width=0.7\columnwidth]{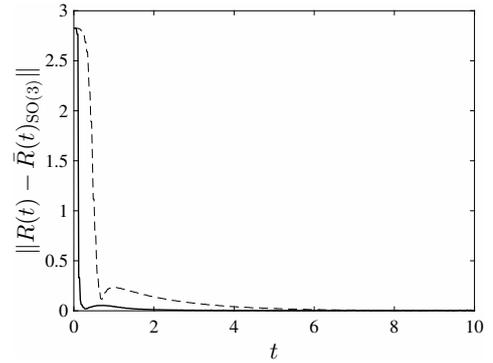}}
    }
    \centerline{
        \subfigure[Gyro bias estimation error]{\includegraphics[width=0.7\columnwidth]{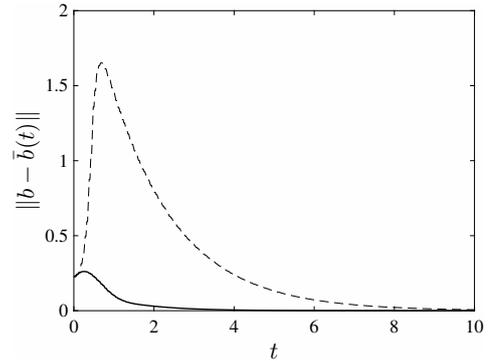}}
    }
\caption{Experimental results: the attitude estimation error $\|R(t)-\bar R(t)_{\SO}\|$, after polar-decomposition, and the gyro bias estimation error $\|b-\bar b(t)\|$ are compared between the proposed observer \eqref{obs:quad:var} (thick,solid) with the Mahony observer (thin,dashed)}\label{fig:exp}
\end{figure}

 For the sake of convenience, let $\bar R (t) \in \mathbb R^{3\times 3}$ and $\bar R_{\SO}(t) \in \SO$ denote the attitude estimate trajectory generated by our observer and its $\SO$ factor obtained through polar decomposition, respectively.  The initial estimates are given by $\bar R(0) = R(0) \exp(0.99\pi \hat e_3)$ and $\bar b(0)=(0,0,0)$,  where $R(t)$ denotes the  trajectory of true orientation. 
The observer gains are selected as $k_P=2.5$, $k_I=1.5$.
The corresponding results are illustrated in Figure \ref{fig:exp}, where the attitude estimation error (after polar decomposition), and the gyro bias estimation error are plotted. 
These experimental results exhibit qualitative behaviors consistent with the numerical simulations. 
For the selected initial condition near the boundary of the region of attraction of the Mahony observer, the Mahony observer with the same observer gains yields a slower initial convergence in the attitude estimate and a larger overshoot in the gyro bias estimate whereas the proposed observer shows more desirable convergence rates.

\section{CONCLUSION}
We have successfully designed a 12-dimensional, global, continuous, and exponentially convergent observer for attitude and gyro bias of a rigid body. 
This observer overcomes the topological restriction on $\SO$ completely by constructing it in the embedding space, and it eliminates fundamental drawbacks of other geometry-free attitude observers.  
The desirable properties of the proposed observers are illustrated by experimental results based on visual landmarks. 
Future works include extending this results to the special Euclidean group for concurrent estimation of the position and the attitude of a rigid body, and integrating this observer with an attitude controller to show stability of the combined system.  We plan to apply the result to drone control (\cite{ChEu17}) and to combine with deep learning (\cite{CaCh18}).


\bibliography{wroco2019}

\begin{thebibliography}{18}
\providecommand{\natexlab}[1]{#1}
\providecommand{\url}[1]{\texttt{#1}}
\providecommand{\urlprefix}{URL }
\expandafter\ifx\csname urlstyle\endcsname\relax
  \providecommand{\doi}[1]{doi:\discretionary{}{}{}#1}\else
  \providecommand{\doi}{doi:\discretionary{}{}{}\begingroup
  \urlstyle{rm}\Url}\fi

\bibitem[{Batista et~al.(2012{\natexlab{a}})Batista, Silvestre, and
  Oliveira}]{BaSiOl12b}
Batista, P., Silvestre, C., and Oliveira, P. (2012{\natexlab{a}}).
\newblock A {GES} attitude observer with single vector observations.
\newblock \emph{Automatica}, 48, 388--395.

\bibitem[{Batista et~al.(2012{\natexlab{b}})Batista, Silvestre, and
  Oliveira}]{BaSiOl12c}
Batista, P., Silvestre, C., and Oliveira, P. (2012{\natexlab{b}}).
\newblock Globally exponentially stable cascade observers for attitude
  estimation.
\newblock \emph{Control Engineering Practice}, 20, 148--155.

\bibitem[{Batista et~al.(2012{\natexlab{c}})Batista, Silvestre, and
  Oliveira}]{BaSiOl12a}
Batista, P., Silvestre, C., and Oliveira, P. (2012{\natexlab{c}}).
\newblock Sensor-based globally asymptotically stable filters for attitude
  estimation: analysis, design, and performance evaluation.
\newblock \emph{IEEE Trans. Automatic Control.}, 57, 2095--2100.

\bibitem[{Berkane et~al.(2017)Berkane, Abdessameud, and Tayebi}]{BeAbTa17}
Berkane, S., Abdessameud, A., and Tayebi, A. (2017).
\newblock Hybrid attitude and gyro-bias observer on {SO}(3).
\newblock \emph{IEEE Trans. Automatic Control}, 62, 6044--6050.

\bibitem[{Caterini and Chang(2018)}]{CaCh18}
Caterini, A. and Chang, D. (2018).
\newblock \emph{Deep Neural Networks in a Mathematical Framework}.
\newblock Springer.

\bibitem[{Chang(2011)}]{Ch11}
Chang, D. (2011).
\newblock A simple proof of the {P}ontryagin maximum principle on manifolds.
\newblock \emph{Automatica}, 47, 630--633.

\bibitem[{Chang(2018{\natexlab{a}})}]{Ch18b}
Chang, D. (2018{\natexlab{a}}).
\newblock Observer-based controller design for systems on manifolds in
  {E}uclidean space.
\newblock \emph{Proc. 2018 57th Annual Conference of the Society of Instrument
  and Control Engineers of Japan (SICE)}, 573--578.

\bibitem[{Chang(2018{\natexlab{b}})}]{Ch18}
Chang, D. (2018{\natexlab{b}}).
\newblock On controller design for systems on manifolds in {E}uclidean space.
\newblock \emph{Int J Robust Nonlinear Control}, 28, 4981--4998.

\bibitem[{Chang and Eun(2017)}]{ChEu17}
Chang, D. and Eun, Y. (2017).
\newblock Global chartwise feedback linearization of the quadcopter with a
  thrust positivity preserving dynamic extension.
\newblock \emph{IEEE Trans. Automatic Control}, 62, 4747 -- 4752.

\bibitem[{Faraut(2008)}]{Fa08}
Faraut, J. (2008).
\newblock \emph{Analysis on Lie Groups}.
\newblock Cambridge University Press, New York.

\bibitem[{Garrido-Jurado et~al.(2014)Garrido-Jurado, {n}oz Salinas,
  Madrid-Cuevas, and Mar\'{i}n-Jim\'{e}nez}]{GarMun14}
Garrido-Jurado, S., {n}oz Salinas, R.M., Madrid-Cuevas, F.J., and
  Mar\'{i}n-Jim\'{e}nez, M.J. (2014).
\newblock Automatic generation and detection of highly reliable fiducial
  markers under occlusion.
\newblock \emph{Pattern Recogn.}, 47, 2280--2292.

\bibitem[{Grip et~al.(2012)Grip, Fossen, Johansen, and Saberi}]{GrFoJoSa12}
Grip, H., Fossen, T., Johansen, T., and Saberi, A. (2012).
\newblock Attitude estimation using biased gyro and vector measurements with
  time-varying reference vectors.
\newblock \emph{IEEE Trans. Automatic Control}, 57, 1332--1338.

\bibitem[{Grip et~al.(2015)Grip, Fossen, Johansen, and Saberi}]{GrFoJoSa15}
Grip, H., Fossen, T., Johansen, T., and Saberi, A. (2015).
\newblock Globally exponentially stable attitude and gyro bias estimation with
  application to {GNSS}/{INS} integration.
\newblock \emph{Automatica}, 51, 158 --166.

\bibitem[{Mahony et~al.(2008)Mahony, Hamel, and Pflimlin}]{Ma08}
Mahony, R., Hamel, T., and Pflimlin, J.M. (2008).
\newblock Nonlinear complementary filters on the special orthogonal group.
\newblock \emph{IEEE Trans. Automatic Control}, 53, 1203--1218.

\bibitem[{Martin and Sala\"{u}n(2007)}]{MaSa07}
Martin, P. and Sala\"{u}n, E. (2007).
\newblock Invariant observers for attitude and heading estimation from low-cost
  inertial and magnetic sensors.
\newblock \emph{Proc. IEEE Conference on Decision and Control}, 1039--1045.

\bibitem[{Martin and Sala\"{u}n(2010)}]{MaSa10}
Martin, P. and Sala\"{u}n, E. (2010).
\newblock Design and implementation of a low-cost observer-based attitude and
  heading reference system.
\newblock \emph{Control Engineering Practice}, 18, 712--722.

\bibitem[{Martin and Sarras(2017)}]{MaSa17}
Martin, P. and Sarras, I. (2017).
\newblock A global observer for attitude and gyro biases from vector
  measurements.
\newblock \emph{IFAC PapersOnLine}, 50-1, 15409--15415.

\bibitem[{Wu et~al.(2015)Wu, Kaufman, and Lee}]{WuKaLe15}
Wu, T.H., Kaufman, W., and Lee, T. (2015).
\newblock Globally asymptotically stable attitude observer on {SO}(3).
\newblock \emph{Proc. IEEE American Control Conference}, 2165--2168.

\end{thebibliography}

\appendix
\section{Mathematical preliminaries}
This appendix contains mathematical preliminaries to help the reader understand the main results of the paper. The usual Euclidean inner product is exclusively used for vectors and matrices in this paper, i.e.  $\langle A, B \rangle = \sum_{i,j}A_{ij}B_{ij} = \operatorname{tr}(A^TB)$ for any two matrices  of equal size.  The norm induced from this inner product, which is called the Frobenius or Euclidean norm, is exclusively used for vectors and matrices. Let $\Sym$ and $\Skew$ denote the symmetrization operator and the skew-symmetrization operator, respectively, on square matrices, which are defined by
\[
\Sym (A) = \frac{1}{2}(A+A^T), \quad \Skew (A) = \frac{1}{2}(A-A^T)
\]
for any square matrix $A$.
Then,
\[
A  = \Sym (A) + \Skew(A), \quad \langle \Sym (A), \Skew (A) \rangle = 0.
\]
 Namely, 
\[
\mathbb R^{n\times n} = \Sym (\mathbb R^{n\times n} ) \oplus \Skew (\mathbb R^{n\times n} )
\]
with respect to the Euclidean inner product. 
Let $\SO$ denote the set of all $3\times 3$ rotation matrices, which is defined as $\SO = \{ R\in \mathbb R^{3\times 3} \mid R^T R =I, \det R=1\}$.  Let $\mathfrak{so}(3)$ denote the set of all $3\times 3$ skew symmetric matrices, which is defined as $\mathfrak{so}(3) = \{ A \in \mathbb R^{3\times 3} \mid A^T+ A = 0 \}$. The hat map $\wedge : \mathbb R^3 \rightarrow \mathfrak{so}(3)$ is defined by
\begin{equation}\label{def:Omega}
\hat \Omega = \begin{bmatrix}
0 & -\Omega_3 &\Omega_2 \\
\Omega_3 & 0 & -\Omega_1\\
-\Omega_2 & \Omega_1 & 0
\end{bmatrix}
\end{equation}
for $\Omega = (\Omega_1, \Omega_2,\Omega_3) \in \mathbb R^3$.  The inverse map of the hat map is called the vee map and denoted  $\vee$ such that $(\hat \Omega)^\vee = \Omega$ for all $\Omega \in \mathbb R^3$ and $(A^\vee)^\wedge = A$ for all $A\in \mathfrak{so}(3)$.  
\begin{lemma}\label{lemma:prelim}
1. $\langle R A, R B\rangle  = \langle A, B\rangle = \langle AR, BR\rangle$   for all $R \in \SO$ and $A,B \in \mathbb R^{3\times 3}$.

2. $\lambda_{\rm min}(A^TA)\|B\|^2 \leq \langle AB,AB\rangle \leq \lambda_{\rm max}(A^TA)\|B\|^2$ for all $A \in \mathbb R^{n\times m}$ and $B\in \mathbb R^{m \times \ell}$.

3. $\langle \hat x, \hat y \rangle = 2 \langle x, y\rangle$ for all $x, y \in \mathbb R^3$.

4. $\|A\|^2 = \|\Sym(A)\|^2 + \|\Skew(A)\|^2$ for all $A\in \mathbb R^{n\times n}$.

5. $\|AB\| \leq \|A\| \|B\|$ for all $A \in \mathbb R^{n\times m}$ and $B\in \mathbb R^{m \times \ell}$.

6. $x\times y = (y x^T - x y^T)^\vee $ for all $x,y\in \mathbb R^3$.

7. $\max_{R_1, R_2 \in \SO}\| R_1 - R_2\| = 2\sqrt 2$. 
\end{lemma}

\end{document}